\documentclass[10pt,fleqn,reqno,twoside]{amsart}

\usepackage[mathscr]{eucal}
\usepackage{amsmath}
\usepackage[body={12.54125cm,19.75cm},truedimen]{geometry}

\makeatletter
\def\article@logo{%
  \set@logo{%
   Hyper{\TeX}ed arXival version prepared and 
   updated$^{\protect\hyperlink{update}{(*)}}$ February 2002.\\
   Originally published in Comment.\ Math.\ Helvetici \textbf{62} %
   (1987) 630--645.
  }%
}
\def\eonly@logo{%
  \set@logo{%
   Hyper{\TeX}ed arXival version prepared and 
   updated$^{\protect\hyperlink{update}{(*)}}$ February 2002.\\
   Originally published in Comment.\ Math.\ Helvetici \textbf{62} %
   (1987) 630--645.
  }%
}
\def\publname{Hyper{\TeX}ed arXival version prepared and 
   updated$^{\protect\hyperlink{update}{(*)}}$ February 2002.}
\def\pageinfo{\relax}
\def\volinfo{\relax}
\def\ISSN{\relax}
\renewcommand\PII[1]{\def\@PII{#1}}
\PII{Originally published in Comment.\ Math.\ Helvetici \textbf{62} %
   (1987) 630--645.}
\def\@setcopyright{\relax}
\def\@settitle{%
\noindent
  \baselineskip12\p@\relax
    \bfseries
\large
  \@title
}
\def\@setauthors{%
  \begingroup
\mdseries\scshape
  \trivlist
\footnotesize \@topsep30\p@\relax
  \advance\@topsep by -\baselineskip
  \item\relax
  \andify\authors
  {\normalsize\authors}%
  \endtrivlist
  \endgroup
}
\def\@secnumfont{\bfseries}
\def\section{\@startsection{section}{1}%
  \z@{1.2\linespacing\@plus\linespacing}{.7\linespacing}%
  {\normalfont\bfseries\noindent\S}}
\def\currentvolume{\relax}
\def\@adjustvertspacing{%
  \bigskipamount1.0\baselineskip 
  \medskipamount\bigskipamount \divide\medskipamount\tw@
  \smallskipamount\medskipamount \divide\smallskipamount\tw@
  \abovedisplayskip\bigskipamount  
  \belowdisplayskip \abovedisplayskip
  \abovedisplayshortskip\abovedisplayskip
  \advance\abovedisplayshortskip-1\abovedisplayskip
  \belowdisplayshortskip\abovedisplayshortskip
  \advance\belowdisplayshortskip 1\smallskipamount
  \jot\baselineskip \divide\jot 4 \relax
}
\renewenvironment{proof}[1][\quad Proof]{%
\normalfont \topsep6\p@\@plus6\p@\relax
\trivlist
\item[\hskip\labelsep
      \itshape
  #1\@addpunct{.}]\ignorespaces
}{\endtrivlist\@endpefalse
}
\makeatother

\usepackage{hyperref}

\newtheoremstyle{cmhthm}%
{\baselineskip plus.2\baselineskip minus.2\baselineskip}{12pt}{\itshape}%
{\parindent}{\mdseries}{.}{ }%
{}

\theoremstyle{cmhthm}
\newtheorem{theorem}{THEOREM}[section]
\newtheorem{corollary}[theorem]{COROLLARY}
\newtheorem{proposition}[theorem]{PROPOSITION}
\newtheorem{lemma}[theorem]{LEMMA}
\newtheorem{definition}[theorem]{DEFINITION}
\newtheorem*{theorem*}{THEOREM}
\newtheorem*{theorem 2.3}{\protect\hyperlink{Theorem 2.3}{THEOREM 2.3}}
\newtheorem*{theorem 2.5}{\protect\hyperlink{Theorem 2.5}{THEOREM 2.5}}
\newtheorem*{corollary 2.6}{\protect\hyperlink{Corollary 2.6}{COROLLARY 2.6}}

\newtheoremstyle{cmhrmk}%
{\baselineskip plus.2\baselineskip minus.2\baselineskip}{12pt}{\upshape}%
{\parindent}{\itshape}{.}{ }%
{}
\theoremstyle{cmhrmk}
\newtheorem{remark}[theorem]{Remark}
\newtheorem{remarks}[theorem]{Remarks}

\newtheoremstyle{cmhconst}%
{\baselineskip plus.2\baselineskip minus.2\baselineskip}{12pt}{\upshape}%
{\parindent}{\upshape}{.}{ }%
{}
\theoremstyle{cmhconst}%
\newtheorem{hypothesis}[theorem]{HYPOTHESIS}
\newtheorem{construction}[theorem]{CONSTRUCTION}
\newtheorem{example}[theorem]{EXAMPLE}
\newtheorem{machinery}[theorem]{MACHINERY}

\newtheoremstyle{addenda}%
{\baselineskip plus.2\baselineskip minus.2\baselineskip}{12pt}{\upshape}%
{}{\upshape}{}{ }%
{\thmname{#1}\thmnumber{}\thmnote{}}%
\theoremstyle{addenda}%
\newtheorem{addenda}{$^{(*)}$ADDENDA (February 2002).}

\def\Ks{\mathscr K}

\def\Ps{\mathscr P}

\def\Us{\mathscr U}

\def\C{\mathbf{C}}
\def\H{\mathbf{H}}

\def\R{\mathbf{R}}
\def\Z{\mathbf{Z}}

\def\iq{\mathbf{i}}
\def\jq{\mathbf{j}}
\def\kq{\mathbf{k}}
\def\pq{\mathbf{p}}
\def\qq{\mathbf{q}}
\def\uq{\mathbf{u}}

\def\xq{\mathbf{x}}
\def\yq{\mathbf{y}}
\def\zero{\mathbf{0}}

\def\cone{\operatorname{cone}}
\def\conj{\operatorname{conj}}
\def\crit{\operatorname{crit}}
\def\GCD{\operatorname{GCD}}
\def\Hopf{\operatorname{Hopf}}
\def\Int{\operatorname{Int}}
\def\lk{\operatorname{\mathit{lk}}}
\def\pos{\operatorname{pos}}
\def\POS{\operatorname{POS}}
\def\neg{\operatorname{neg}}
\def\NEG{\operatorname{NEG}}
\def\Re{\operatorname{Re}}
\def\Im{\operatorname{Im}}
\def\Rev{\operatorname{Rev}}

\def\after{\raise.25ex\hbox{$\scriptstyle\thinspace\circ\thinspace$}}
\def\Bd{\partial}
\def\e{\varepsilon}
\def\from{\colon\thinspace}
\def\phi{\varphi}
\def\s{\sigma}
\def\sub{\subset}

\newcommand{\Span}[1]{\langle\protect #1\rangle}

\begin{document}
\renewcommand{\thefootnote}{\fnsymbol{footnote}}
\setcounter{section}{-1}
\setcounter{addenda}{0}
\mathindent12pt

\title[Splitting the Milnor number]
{Isolated critical points of mappings from 
$\R^4$ to $\R^2$ and a natural\\	
splitting of the Milnor number of a classical fibered link.\\
Part I: Basic theory; examples}

\author{Lee Rudolph$^{(1)}$}
\thanks{
$\mathstrut^{\protect\hypertarget{update}{(*)}}$\thinspace%
See Addenda~\protect\hyperlink{projected}{1}-%
\protect\hyperlink{more N-R}{4} at the end of the paper.\hfill\break%
\protect\phantom{$2$-p}%
$\mathstrut^{(1)}$\thinspace%
Research partially supported by the Fonds National Suisse.
}
\date{}

\renewcommand{\footnoterule}%
{\vskip5pt \kern -3pt \hrule height 0.4pt width \columnwidth \kern 2.6pt}

\setcounter{page}{630}

\maketitle
\baselineskip=1.05\baselineskip
\markboth{\scshape{lee rudolph}}{Splitting the Milnor number}

\section{Introduction; statement of results}
\label{section 0}

From a fibered link $\Ks = (S^3, K)$ may be 
constructed a field $S_\Ks$ of (not\break%
everywhere tangent) $2$-planes on $S^3$. When $\Ks = \Ks_F$
is the link of an isolated
critical point of a map $F\from  \R^4\to\R^2$, $S_\Ks$ is 
essentially the field of kernels of $DF$.
Homotopically, $S_\Ks$ determines integers $\lambda(\Ks)$ and 
$\rho(\Ks)$.

\begin{theorem*}$\lambda(\Ks) + \rho(\Ks) = \mu(\Ks)$.
\end{theorem*}

Here $\mu(\Ks)$ is the Milnor number of $\Ks$, that is, 
the rank of the first homology\break%
of the fiber surface of $\Ks$. At least in some 
cases, this splitting of $\mu(\Ks)$\break%
corresponds to a geometrically 
natural direct sum decomposition of this homology
group, into subgroups corresponding to ``negative'' 
(or ``left-handed'') and\break%
``positive'' (or ``right-handed'') parts of the fiber 
surface, \cite{Rudolph:isocp2}. For instance, if $\Ks$\break%
is a closed positive braid (such as the Lorenz 
links of dynamical systems \cite{Birman-Williams}, or\break%
the overlapping class of links of complex plane 
curve singularities -- links $\Ks_F$\break%
where $F\from\C^2\to \C$ is complex-analytic), then 
$\lambda(\Ks) = 0$ and $\rho(\Ks) = \mu(\Ks)$,
substantiating the intuition that such links are 
as positive as they can be. If $\Ks$ is\break%
the figure-$8$ knot, then $\lambda(\Ks) = 1 = \rho(\Ks)$.

The Euler characteristic $1-\mu(\Ks)$ of a fiber 
surface of $\Ks$ can be computed by
correctly counting the singularities of a 
vectorfield on the surface. There is a sense
in which the extra information in the splitting 
$\mu(\Ks) = \lambda(\Ks)+ \rho(\Ks)$ comes from
making this calculation ``all around the circle'' 
of fiber surfaces.

This is the first of several papers devoted to 
the study of $\lambda$, $\rho$, and related
invariants. In this paper I develop the basic 
theory, and compute a number of
examples.

More specifically, in \S\ref{section 1}, I construct the field 
$S_\Ks$, and two related (tangent)
$2$-plane fields $T^{\pm}_\Ks$, from an open-book structure on 
$S^3$ of type $\Ks$. A standard\break%
\newpage\noindent
parallelization of $S^3$ permits one to consider 
$S_\Ks$ as a map $S^3 \to S^2\times  S^2$ and $T^{\pm}_\Ks$ as\break%
maps $S^3\to S^2$; the Hopf invariants 
$(-\lambda(\Ks), \rho(\Ks))$
and $\tau^{\pm}(\Ks)$ of these maps
depend only on $\Ks$. In \S\ref{section 2}, the following results 
are obtained.

\begin{theorem 2.3} For any fibered link $\Ks$, 
$\tau^+(\Ks) - \tau^-(\Ks) = 1 - \mu(\Ks)$.
\end{theorem 2.3}

\par

\begin{theorem 2.5} For any fibered link $\Ks$,
$\lambda(\Ks) = -\tau^+(\Ks)$ and $\rho(\Ks) =\break%
\tau^-(\Ks) + 1$.
\end{theorem 2.5}

\par

\begin{corollary 2.6} For any fibered link $\Ks$, 
$\lambda(\Ks) + \rho(\Ks) = \mu(\Ks)$.
\end{corollary 2.6}

In \S3, 
I exploit the close relationship between 
fibered links and isolated\break%
critical points of mappings from $\R^4$ to $\R^2$ to give 
another way to calculate $\lambda$ and $\rho$\break%
(in theory), which is used in examples in \S4. 

In \cite{Rudolph:isocp2}, 
$\lambda$ and $\rho$ are computed for closed strict 
generalized homogeneous\break 
braids. 
In \cite{Rudolph:mbob}, the pair $\{\tau^{\pm}(\Ks)\}$ 
of invariants of a fibered link $\Ks$ in $S^3$ is\break%
generalized to a set $\{\tau^s(\Ks):s\in \{+,-\}^m\}$
of invariants of a fibered link $\Ks$ of $m$\break%
components in any closed $3$-manifold. Work with 
Walter Neumann (\cite{Neumann-Rudolph:unfoldings}, 
\cite{Neumann-Rudolph:computing}) extends the definition 
of $\lambda$ and $\rho$ to higher dimensions, and studies 
their behavior (in dimension $3$) under various geometric 
operations (cabling, connected sum, Murasugi sum, Stallings 
twists); it is shown, in particular, that \emph{for many 
pairs of\break%
 fibered links $\Ks, \Ks'$ in $S^3$, 
$\{\lambda(\Ks), \rho(\Ks)\}\ne \{\lambda(\Ks'), \rho(\Ks')\}$ 
although $\Ks$ and $\Ks'$\break%
have identical Seifert forms and algebraic 
monodromy}. A \hypertarget{projected-backlink}{projected 
future paper}\hyperlink{projected}{$^{(*)}$}
generalizes $\lambda$ and $\rho$ to non-fibered links, 
using \cite{Gabai:MSRI}.

\hypertarget{regret-backlink}{An interesting phenomenon} 
(still unexplained as of December, 1986) 
is that in no known example is $\lambda(\Ks)$ or $\rho(\Ks)$ 
negative.\hyperlink{regret}{$^{(*)}$}

\bigskip
\bigskip

\section{\protect\label{section 1}%
Some plane fields associated to a classical 
fibered link; the invariants\break%
\protect\phantom{.}\hskip-1.5em 
$\lambda,\rho, \text{ and }\tau^{\pm}$}

\medskip

The $3$-sphere $S^3$ is the boundary of the unit $4$-disk 
$D^4$ of $\C^2$. The $1$-sphere $S^1$ is\break%
the boundary of the unit disk $D^2$ of $\C$. These 
spaces are oriented by the usual
conventions; $\R^4$ (resp.\ $\R^2$) is the oriented real 
vectorspace underlying $\C^2$ (resp.\ $\C$); 
$\H$, the quaternions, is oriented by the usual 
identification of its underlying
vectorspace with $\R^4$. 
The oriented span of a $k$-frame $(U_1,\dots, U_k)$ is\break%
$\Span{U_1,\dots,U_k}$.

A \emph{link} $\Ks$ is a pair $(S^3,K)$ where $K$ is an 
oriented closed smooth $1$-submanifold
of $S^3$ (if $K$ is connected, $\Ks$ is also called a 
\emph{knot}); $\Ks$ is \emph{fibered} if there exist a\break%
\newpage\noindent
closed disk-bundle neighborhood $N(K)$ of $K$ in $S^3$ 
with a smooth trivialization
\[
\psi\from N(K)\to D^2
\]
and a smooth fibration over $S^1$ of the link 
exterior $E(K) = S^3\setminus \Int N(K)$
\[
\phi\from  E(K)\to S^1
\]
such that $\psi\,|\,\Bd N(K)=\phi\,|\,\Bd N(K)$.  Such maps
$\psi,\phi$ glue together to give a\break%
piecewise smooth map
\[
\pi\from  S^3\to D^2
\]
which is called an \emph{open-book structure on $S^3$}, of 
\emph{type $\Ks$}. (This summary treatment
follows \cite{Neumann-Rudolph:unfoldings}; 
cf.\ also \cite{Kauffman-Neumann}.) If $\Ks$ is a fibered 
link, then all open-book structures of
type $\Ks$ are equivalent in quite a strong sense 
\cite{Kauffman-Neumann}, and the ambient isotopy type of\break%
$\Ks$ determines the ambient isotopy type of any 
\emph{fiber surface} of $\Ks$ -- i.e., a fiber of
(any) $\phi$.  The \emph{Milnor number} $\mu(\Ks)$
of $\Ks$ is the rank of the first homology of a\break%
fiber surface of $\Ks$.

Associated to a fibered link $\Ks$ are maps $S_\Ks$, 
$T^{+}_\Ks$, and $T^{-}_\Ks$ from $S^3$ to\break%
$G = G^{+}(2,4)$, the Grassmann manifold of oriented 
$2$-planes in $\R^4$, which will be\break%
defined in terms of certain vectorfields.

Let $\pi$ be an open-book structure of type $\Ks$. Then 
$\pi/|\pi|\from  S^3\setminus K \to S^1$ is a\break%
fibration (by fiber surfaces to which open 
collars are attached piecewise-\break%
smoothly) extending 
$\phi$. It is always possible to 
take $\pi$ such that $\pi/|\pi|$ is smooth,
and we do so. Let $U$ be the field of unit tangent 
vectors on $S^3 \setminus K$ perpendicular to\break%
the plane $\ker(D(\pi/|\pi|))$, so oriented as to point 
to the positive side of the fiber
surfaces. Let $V$ be the field of unit tangent 
vectors on $N(K)$ which span the line\break%
$\ker (D\pi)$, so oriented as to induce the given 
orientation on $K = \pi^{-1}(0)$. In\break%
$N(K) \setminus K$, where both are defined, the 
vectorfields $U$ and $V$ are orthogonal; let\break%
$U\times V$ be the field of unit vectors, tangent to 
$S^3$, such that the orthonormal\break%
$3$-frame $(U, V, U\times  V)$ gives the orientation of 
$S^3$. Finally, let $W$ be the field of
unit inward normal vectors on $S^3$ (so the 
orthonormal $4$-frame $(U, V, U\times V, W)$
gives the orientation of $\C^2$).

\begin{definition} For $Q \in E(K)$, 
$S_\Ks(Q) = T^{+}_\Ks(Q) = T^{-}_\Ks(Q)$
is the oriented tangent plane to the fiber surface 
$\phi^{-1}(\phi(Q))$ through $Q$. For $Q \in N(K) \setminus K$,
$S_\Ks(Q) = 
\Span{V(Q), |\pi(Q)|(U \times V)(Q) + (1 - |\pi(Q)^2)^{1/2}W(Q)}$
and $T^{\pm}_\Ks(Q) =\break%
\Span{|\pi(Q)|V(Q)\mp (1 - |\pi(Q)^2)^{1/2}U(Q),(U \times V)(Q))}$. 
For $Q \in K$, 
$S_\Ks(Q) =\break%
\Span{V(Q), W(Q)}$ and 
$T^{\pm}_\Ks(Q)$ is the plane tangent to $S^3$ 
and orthogonal to $V$ at $Q$, so\break%
oriented that $\pm V$ points to its positive side.
\end{definition}

\newpage
It is easy to see that the maps $S_\Ks$ and $T^{\pm}_\Ks$ are 
continuous, and that their\break%
homotopy classes in $\pi_3(G)$ depend only on $\Ks$ and 
not on the choice of $\pi$ (so the\break%
notation is not too abusive).

It is well known that $G$ is diffeomorphic to $S^2\times S^2$, 
so that $\pi_3(G)$ is\break%
isomorphic to $\Z\oplus\Z$. 
To obtain an explicit isomorphism (and, thus, integer 
invariants of $\Ks$), we use quaternions to give an 
explicit diffeomorphism, with pleasant properties 
which will facilitate later computations. Let 
$1, \iq, \jq, \kq$ be the standard orthonormal 
basis of $\H$; in terms of the 
identification of $\H$ with $\C^2$
already established, we have $1 = (1,0)$, 
$\iq = (i,0)$, $\jq = (0,1)$, and $\kq = (0, i)$ 
(note that we distinguish the quaternion $\iq$ from the 
complex number $i$ by boldface). The
tangent space of $S^3$ at $1$ is $\Span{\iq, \jq, \kq}$, 
the \emph{pure quaternions}; 
let $S^2 = S^3\cap \Span{\iq, \jq, \kq}$ be
the oriented sphere of unit pure quaternions. 
Let $\Ps\from \H\to\Span{\iq,\jq,\kq}$ be the ``pure
part'' mapping; write $\conj \from  H\to \H$ for quaternionic 
conjugation; let 
$\Us\from  H\setminus\{0 \} \to S^3:Q\mapsto Q/\|Q\|$. For 
$Q\in\H$, let $L_Q\from\H\to\H$ (resp., $R_Q$) be 
the real-linear map $A\mapsto QA$ (resp.\ $A\mapsto AQ$).

It is a fact that $S^2$ consists precisely of the 
square roots of $-1\,\in\,\H$. Hence each\break%
$\pq \in S^2$ determines two complex structures on $\R^4$, 
with structure maps respectively
$L_\pq$ and $R_\pq$; call $(\R^4, L_\pq)$ 
(resp.\ $(\R^4, R_\pq)$ the \emph{$\pq$-left} 
(resp.\ \emph{$\pq$-right}) \emph{complex\break%
structure}. 
(The $\iq$-left complex structure is the 
original structure of $\R^4$ as $\C^2$; the\break%
$\iq$-right complex structure is as it were the 
direct sum $\C\oplus \bar\C$.) Each $\pq$-left (resp.
$\pq$-right) complex structure determines a subset of 
$G$, namely, the oriented real\break%
 $2$-planes which are 
left (resp.\ right) $\pq$-stable -- that is, which 
are complex lines in that structure. (Note that if 
a plane is left or right $\pq$-stable with one orientation
then the same plane is left or right $-\pq$-stable 
with the opposite orientation.)\break%
There are, in fact, well-defined maps $l$ and $r$ 
from $G$ to $S^2$ such that, for every $W \in G$, 
$W$ is a complex line in the $l(W)$-left 
complex structure and in the\break%
$r(W)$-right complex structure, 
and in those structures only.

\begin{lemma}\label{lemma 1.2}
If $(A, B)$ is a $2$-frame, then 
$l(\Span{A,B }) = \Us\Ps(B \conj (A))$ and
$r(\Span{A,B }) = \Us\Ps(\conj(A)B)$. 
In complex coordinates, if $A = (z_1, w_1), B =\break%
(z_2, w_2)$, then

\begin{flalign*}
l(\langle{A,B}\rangle)=\Us\Ps(\bar z_1 z_2+\bar w_1 w_2, z_1 w_2-w_1 z_2),\\
r(\langle{A,B}\rangle)=\Us\Ps(\bar z_1 z_2+ w_1\bar w_2, %
z_1\bar w_2-\bar w_1 z_2).
\end{flalign*}

\noindent
The pair $(l,r)\from G\to S^2 \times S^2$ is a 
diffeomorphism.
\end{lemma}

The compositions $(l, r)\after S_\Ks$, 
$(l, r)\after T^{+}_\Ks$, and $(l, r)\after T^{-}_\Ks$
map $S^3$ to $S^2\times S^2$ and\break%
so provide elements of $\pi_3(S^2 x S^2) = \pi_3(S^2)\oplus\pi_3(S^2)$. 
Recall that the \emph{Hopf}\break%
\newpage\noindent
\emph{invariant} $H(g)$ of a continuous 
map $g:S^3\to S^2$ can 
be defined as follows: let $\pq$ and\break%
$\qq$ be 
distinct regular values of $g_1$, a map 
homotopic to $g$ which is smooth near
$g_1^{-1}(\pq, \qq)$; then $H(g)$ is the linking number 
of the smooth links $g_1^{-1}(\pq)$ and\break%
$g_1^{-1}(\qq)$, where $g_1^{-1}(\pq)$ 
is oriented so that, if $D$ 
is a small oriented normal disk
intersecting it once positively, then 
$g_1 | D \from D \to S^2$ preserves orientation 
(and similarly for $g_1^{-1}(\qq)$). The Hopf invariant of
\[
\Hopf\from S^3\to S^2:%
(z, w)\mapsto ([|z|^2-|w|^2]i, 2iz\bar w),
\]
the \emph{complex} (as opposed to conjugate-complex) 
\emph{Hopf fibration}, is $+1$. The Hopf
invariant of maps induces an isomorphism 
(denoted by the same name and letter)
$H\negthinspace:\negthinspace\pi_3(S^2)\to \Z$. 
If ${-}\negthinspace: S^2\to S^2$ is the antipodal 
map, then $H(-\after g) = H(g)$ for any\break%
$g\from  S^3\to S^2$, since fibers of $-\after g$ 
are fibers of $g$ with orientation reversed, and
linking number is bilinear.

\begin{lemma}\label{lemma 1.3}
If $T: S^3\to G$ is such that $T(Q)$ is
tangent to $S^3$ at $Q$ for every
$Q \in S^3$, then $H(r\after T) = 1 + H(l\after T)$.
\end{lemma}

\smallskip

\begin{proof} We first check a special case in which 
$H(l\after T) = 0$. The tangent\break%
space to $S^3$ at $Q$ is $\Span{\iq Q, \jq Q, \kq Q}$. 
Let $T(Q) = \Span{\jq Q, \kq Q}$; 
then $(l\after T)(Q) = \iq$ is
constant, so $H(l\after T) = 0$. 
On the other hand, $(r\after T)(Q) = \Us\Ps(\conj(\jq Q)\kq Q) =
Q^{-1}(-\jq\kq)Q = -\Hopf (Q)$, 
so $H(r\after T) = 1$.

In general, let $T(Q) = \Span{\pq(Q)Q, \qq(Q)Q}$.  Then

\smallskip

\begin{flalign*}
l(T(Q)) &= \qq(Q)Q\conj(\pq(Q)Q) = -\qq(Q)\pq(Q) = \pq(Q)\qq(Q),\\
r(T(Q)) &= \conj(\pq(Q)Q)\qq(Q)Q = -Q^{-1}\pq(Q)\qq(Q)Q,
\end{flalign*}
so $r\after T = Ad\after (id,-\after T)$, 
where $Ad\from S^3\times S^2\to S^2: (Q, x) \to  Q^{-1}xQ$. Thus
\begin{align*}
H(r\after T) &= H(Ad\after (id, -\after T)) = H(Ad_\#([id],[-\after T])) \\
             &= H(Ad_\#([id], [*])) + %
                H([-\after T])H(Ad_\#([*],[\Hopf ]))\\
             &= 1+H(l\after T)
\end{align*}
by the special case and the sentence preceding 
the lemma.
\end{proof}

\smallskip

\begin{definition}
By $\lambda(\Ks)$, $\rho(\Ks)$, and 
$\tau^{\pm}(\Ks)$ will be denoted the integers
$-H(l\after S_\Ks)$ \emph{(}note the sign\emph{!)}, 
$H(r\after S_\Ks)$, and 
$H(l\after T^{\pm})$, respectively.
\end{definition}

\newpage
\section{Braided open-book structures; relations among 
$\tau^{\pm}$, $\lambda$, $\rho$, and $\mu$.}\label{section 2}

Let $O=\{(z, 0): |z| = 1\} = \Hopf^{-1}(\iq)$, 
$O'=\{(0,w): |w| = 1\} = \Hopf^{-1}(-\iq)$. 
Let $D$ be a round disk on $S^2$ centered at $\iq$, 
$N = N(O) = \Hopf^{-1}(D)$. The map\break%
$S^3\setminus O' \to O: (z, w) \mapsto (z/|z|, O)$ is a fibration 
by open great hemi-$2$-spheres; its\break%
restriction to $N$ presents $N$ as a disk-bundle 
neighborhood of $0$, with fibers\break%
\emph{meridional disks} of $O$. Let $R$ be the oriented unit 
tangent vectorfield to the oriented fibers of Hopf.

Let $\Ks' = (S^3, K')$ be any link. Then if $1 >\e > 0$,
there is an ambient isotopic link $\Ks = (S^3, K)$ such 
that $K$ is contained in $\Int N$ and, at each point of $K$, the\break%
component of $R$ along the positively directed 
tangent line to $K$ is at least $1 - \e$.
(Apply the classical lemma of Alexander to find 
an isotopy carrying $K'$ onto a closed braid with axis $0'$, 
on some number $n>0$ of strings; then use a ``radial''
isotopy in the open solid torus $S^3\setminus 0'$ to make 
this closed braid lie arbitrarily\break%
$C^1$-close to $O$, the core of the solid torus, 
and note that $R|O$ is the field 
of unit tangent vectors to $0$.) In fact, $K$ can be taken to 
have a disk-bundle\break%
neighborhood $N(K)$ that intersects each 
meridional disk of $O$ in a union of $n$ \emph{meridional disks of $K$},
which can be taken to be round (in the spherical geometry
of the meridional disk of $0$), with a trivialization 
$\psi\from N(K) \to  D^2$ such that the component of $R$ 
along the oriented line $\ker(D\psi)$ is at least $1 - \e$.

If, further, $\Ks$ is fibered, then there are such 
a trivialization of $N(K)$ as above, and a fibration 
$\phi$ of $E(K)$ over $S^1$, which glue together to give an 
open-book structure $\pi$ of type $\Ks$ such that $\pi$ is smooth 
and $\iq$ and $-\iq$ are regular values of 
$l\after S$, $r\after S$, $l\after T^{+}$, and 
$r\after T^{-}$ (where $S = S_\Ks$ and 
$T^{\pm}=T^{\pm}_\Ks$ are the plane fields
constructed from $\Ks$ as in \S\ref{section 1}). We will call such 
a $\pi$ a \emph{braided open-book structure}.

\begin{lemma}\label{lemma 2.1}
Let $\pi$ be a braided open-book 
structure. Then: $(1)$~each of\break%
$(l\after S)^{-1}(-\iq)$, 
$(r\after S^{-1}(\iq)$, $(l\after T^{+})^{-1}(\mp\iq)$, and 
$(r\after T^{-})^{-1}(\pm\iq)$ has empty intersection
with $N(K)$ \emph{(}in particular, it is homologous to $0$
in $N(K)$\emph{);} $(2)$~the naturally oriented\break%
$1$-submanifold $N(K)\cap (l\after S)^{-1}(\iq)$
\emph{(}resp.\ $N(K)\cap (r\after S)^{-1}(-\iq)$\emph{;}
$N(K)\cap (l\after T^{\pm})^{-1}(\pm\iq)$\emph{;}
$N(K)\cap (r\after T^{\pm})^{-1}(\mp\iq)$\emph{)}
of $N(K)$ is homologous to $K$ \emph{(}resp.\ $-K$\emph{;} 
$\pm K$; $\mp K$\emph{)} in $N(K)$.
\end{lemma}

\begin{proof} By using $R$ in place of $V$, construct plane 
fields $\tilde S$ and $\tilde T$ on $N(K)$.\break%
Note that $W = \iq R$. Using \ref{lemma 1.2}, one calculates that 
$l\after\tilde S$, $r\after\tilde S$, 
$l\after\tilde T^{\pm}$, and $r\after\tilde T^{\pm}$ have
non-negative $\iq$ component in $N(K)$; it follows 
easily that the $\iq$ components of 
$l\after S$, $r\after S$, 
$l\after T^{\pm}$, and $r\after T^{\pm}$ 
are bounded away from $-1$, establishing the statements in $(1)$.

To verify $(2)$, first note that (the fundamental 
class of) $K$ generates\break%
$H_1(N(K); \Z)$, so what has to be determined 
in each case is what integer multiple\break%
\newpage\noindent
of $K$ the $1$-submanifold in question represents, 
and this is given by its linking
number with the boundary of any one of the 
meridional disks of $K$. An argument
similar to that given for $(1)$ can now be 
applied.
\end{proof}

\begin{definition}\label{def 2.2} 
Let $\pos(\Ks)$ \emph{(}resp.\ $\neg(\Ks)$\emph{)} 
denote the oriented $1$-sub-\break%
manifold $E(K)\cap (l\after S)^{-1}(\iq) = %
E(K)\cap (l\after T^{\pm})^{-1}(\iq)$ 
\emph{(}resp.\ $E(K)\cap (l\after S)^{-1}(-\iq)) = 
E(K)\cap (l\after T^{\pm})^{-1}(-\iq)$\emph{)} of $E(K)$.
\end{definition}

\begin{theorem}\label{theorem 2.3} 
\hypertarget{Theorem 2.3}{}
For any fibered link $\Ks$, 
$\tau^{+}(\Ks) - \tau^{-}(\Ks)=1-\mu(\Ks)$.
\end{theorem}

\begin{proof}
We may assume we have a braided open book 
structure $\pi$ of type\break%
$\Ks$. Then $\tau^{+}(\Ks) = H(l\after T^{+}) = %
\lk((l\after T^{+})^{-1}(\iq), (l\after T^{+})^{-1}(-\iq))$. 
By \ref{lemma 2.1}-2, this equals 
$\lk(K + \pos(\Ks), \neg(\Ks))$. 
Similarly $\tau^{-}(\Ks) = H(l\after T^{-}) = %
\lk((l\after T^{-})^{-1}(\iq),\break%
(l\after T^{-})^{-1}(-\iq))=%
\lk(\pos(\Ks), -K+\neg(\Ks))$. 
Thus $\tau^{+}(\Ks) - \tau^{-}(\Ks) = \lk(K,$\break
$\pos(\Ks)+\neg(\Ks))$. 
Now, the linking number of $K$ with an oriented 1-sub-\break%
manifold $L$ of $S^3 \setminus K$ is equal to the intersection 
number of $L$ with any\break%
Seifert surface of $K$. 
It follows that $\lk(K, \pos(\Ks) + \neg(\Ks))$ 
is the intersection number of $\pos(\Ks)+neg(\Ks)$ 
with a fiber surface $F$ of $\Ks$, that is, the algebraic
number of points of $F$ where the tangent plane to 
$F$ is left $\pm\iq$-stable. These points\break%
are exactly the zeroes of a certain tangent 
vectorfield on $F$ (namely, the\break%
orthogonal projection of $R$ into the tangent plane to $F$), 
and the sum of the\break%
indices of the zeroes of that vectorfield 
equals $1 -\mu(\Ks)$, the Euler characteristic\break%
of $F$. The theorem follows upon observing that the 
multiplicity assigned to a point\break%
of $(\pos(\Ks) + \neg(\Ks))\cap F$ by the orientation of 
$\pos(\Ks) + \neg(\Ks)$ is the index of\break%
that vectorfield at the point.
\end{proof}

\begin{remark}
In some vague sense, the new 
information in the splitting of\break%
$1 -\mu(\Ks)$ as $\tau^{+}(\Ks)-\tau^{-}(\Ks)$ is coming from 
carrying out the vectorfield argument ``all around the circle'' 
of fiber surfaces.
\end{remark}

\smallskip

\begin{theorem}\label{theorem 2.5}\hypertarget{Theorem 2.5}{}
For any fibered link $\Ks$, $\lambda(\Ks) = -\tau^{+}(\Ks)$ 
and $\rho(\Ks) =$\break%
$\tau^{-}(\Ks) + 1$.
\end{theorem}

\smallskip

\begin{corollary}\label{corollary 2.6}\hypertarget{Corollary 2.6}{}
For any fibered link $\Ks$, $\lambda(\Ks)+\rho(\Ks)=\mu(\Ks)$.
\end{corollary}

\begin{proof}[\quad Proof of $\ref{theorem 2.5}$]
By \ref{lemma 2.1}-2, $\lambda(\Ks) = -H(l\after S) = %
-\lk(K + \pos (\Ks), \neg(\Ks)) =\break%
-\tau^{+}(\Ks)$. Just for this proof, let $\POS(\Ks) = 
E(K)\cap (r\after S)^{-1}(\iq)$, 
$\NEG(\Ks) = E(K)\cap (r\after S)^{-1}(-\iq)$. 
Then $\rho(\Ks) = H(r\after S) = \lk(-K + \POS(\Ks), \NEG(\Ks)) = %
H(r \after T^{-}) = H(l\after T^{-}) + 1 = \tau^{-}(\Ks) + 1$ 
(using \ref{lemma 1.3}).
\end{proof}

\section{Isolated critical points; $\lambda$ and $\rho$ as 
intersection numbers}\label{section 3}

Let $U$ be an open neighborhood of a point $\xq \in \R^4$; 
let $f\from U\to \R^2$ be\break%
continuous at $\xq$ and smooth in $U\setminus \{\xq \}$. 
Denote by $Df$ the (real) differential of $f$,\break%
which we take to be a smooth mapping from 
$U\setminus\{\xq \}$ into the space of $2$-by-$4$\break%
matrices, the rows of $Df(\yq )$ being the gradients 
at $\yq $ of the components of $f$. As
usual, $\yq$  is called a \emph{regular point} of $f$ if 
$Df(\yq)$ has rank $2$, a \emph{critical point} otherwise.
Slightly extending the standard usage, we will
call $\xq$ an \emph{isolated critical point} of $f$
if, for $\e>0$ sufficiently small, 
every $\yq$ with $0<\|\xq-\yq\|\le\e$ is a regular point of $f$.
(So, if it happens that $f$ is smooth at $\xq$ and $\xq$ is 
a regular point of $f$, by this usage $\xq$\break%
is also an isolated 
critical point of $f$.)  By ``counting constants'' one finds that the
expected dimension of the set of critical points 
of a smooth mapping is $1$; thus,
the (genuinely critical) isolated critical 
points are unusual, but of correspondingly
great interest.

If $\yq$ is a regular point of $f$, then the matrix 
$Df(\yq)$ considered as an ordered
pair of rows is a $2$-frame, and so $\Span{Df}$ is a 
smooth mapping from the set of\break%
regular points of $f$ into the Grassmann manifold $G$.

\begin{definition}\label{definition 3.1} 
Let $\xq$ be an isolated critical point of $f$. 
For $\e > 0$ small enough that $\yq$ is a regular point 
for $0<\|\xq-\yq\|\le \e$, let 
$E\from S^3\to S^3(\xq,\e):\uq \mapsto \xq + \e\uq$; 
then define $\lambda(f; \xq) = -H(l\after\Span{Df} \after E)$, 
$\rho(f; \xq) = H(r\after\Span{Df} \after E)$.
\emph{(Clearly these do depend only on $f$ and $\xq$.)}
\end{definition}

Two basic facts about isolated critical points 
are relevant here: (A)~when $f$ is\break%
sufficiently well-behaved (e.g., real-polynomial) 
near its isolated critical point $\xq$,
there is an associated ``local link'' $\Ks(f;\xq)$, 
well-defined up to ambient isotopy,
and $\Ks(f;\xq)$ is fibered; (B)~conversely, given a 
fibered link $\Ks in S^3$, there may be
constructed a well-behaved $f_\Ks\from \R^4\to \R^2$ with 
an isolated critical point at $\zero$, such
that $\Ks(f_\Ks; \zero) = \Ks$. More details will be 
recalled shortly.

\begin{remark}\label{remark 3.2}
Milnor \cite[Sect. 10]{Milnor:singular-points} proved (A) for 
polynomial mappings;\break%
Kauffman and Neumann \cite{Kauffman-Neumann} extracted from his proof 
the relevant property of
real-polynomial maps, which they called 
``tameness'', see 3.3. So far as I know,
(B)~was shown first by Looijenga \cite{Looijenga}, and 
rederived in \cite{Kauffman-Neumann} (see 3.7). All this work
is in general dimensions. None of it describes $\lambda$ 
or $\rho$.
\end{remark}

If $\xq$ is a regular point of $f$ and $U$ is a 
sufficiently small open neighborhood of $\xq$,\break%
then $U\cap f^{-1}(f(\xq))$ is a smooth $2$-submanifold of 
$U$. If $\xq$ is an isolated critical
point of (even a smooth) $f$, then this generally 
fails (but not always: cf.\ examples
in \cite{Milnor:singular-points}); all we can say is that, 
for suitably small $U$, the level set $f^{-1}(f(\xq))$\break%
\newpage\noindent
intersects $U$ in a ``2-submanifold with an isolated 
singularity at x''. To give a more
precise description of this singularity, we have 
to impose extra hypotheses on $f$ in
a neighborhood of $\xq$.

\begin{definition}[\cite{Kauffman-Neumann}]%
\label{definition 3.3} Let $\xq$ be an isolated 
critical point of $f$. Then $f$ is tame
at $\xq$ if, for all sufficiently small $\e > 0$,

\hypertarget{hypothesis A}{\emph{(A)}}~the level set 
$f^{-1}(f(\xq))$ intersects $S^3(\xq;\e)$ 
transversely,

\hypertarget{hypothesis B}{\emph{(B)}}~for 
all sufficiently small $\delta=\delta(\e)>0$ the 
intersection $f^{-1}(D^2(f(\xq);\delta))\cap D^4(\xq;\e)$
is a $4$-ball, smooth except for corners along 
$f^{-1}(S^1(f(\xq); \delta))\cap S^3(\xq;\e)$.

\noindent
If $f$ is tame at $\xq$, then \emph{(}for any sufficiently 
small $\e > 0$\emph{)} let $\Ks(f;\xq) =\break%
(S^3,(1/\e)(S^3(\xq;\e) \cap f^{-1}(f(\xq))))$; this is the
\emph{local link} of $f$ at $\xq$.
\end{definition}

\begin{remarks}\label{remarks 3.4}
(1)~As given in \cite{Kauffman-Neumann}, the definition 
of ``tame'' includes the\break%
(inessential) further hypothesis that $f$ is smooth 
at $\xq$. (2)~A mapping can have an\break%
isolated critical point at which it is smooth but 
not tame. (3)~\hyperlink{hypothesis A}{Hypothesis~(A)} of\break%
\ref{definition 3.3} already ensures that $\Ks(f;\xq)$ 
is well-defined (up to ambient isotopy);\break%
\hyperlink{hypothesis B}{Hypothesis~(B)} ensures 
that $\Ks(f; \xq)$ is fibered. 
(4)~As remarked in \cite{Kauffman-Neumann}, the proof
of the ``fibration theorem for real singularities'' 
in \cite{Milnor:singular-points} consists of showing that a
real polynomial mapping is tame at an isolated 
critical point.
\end{remarks}

\begin{construction}\label{construction 3.5}
Let $\pi: S^3\to  D^2$ be a smooth 
open book structure of
type $\Ks$. Define $\cone \pi\from \R^4\to \R^2$
by $(\cone \pi)(\yq) = \|\yq\|\,\,\pi(\yq/\|\yq\|)$ 
if $\yq\ne \zero$, $(\cone\pi)(\zero) = (0,0)$. 
Then the only critical point of $\cone\pi$ is 
$\zero$, $\cone \pi$ is tame at $\zero$, and
$\Ks(\cone\pi;\zero) = \Ks$.
\end{construction}

\begin{remarks}\label{remarks 3.6}
(1)~This construction is a stripped-down version of the original
one in \cite{Looijenga}. (Looijenga showed that, by an 
appropriate choice of $\pi$, $\cone\pi$ can be\break%
taken to be a real polynomial in $\xq$ and $\|\xq\|$, and 
thus real-algebraic, though typically not smooth but merely 
continuous at $\zero$; when, however, $\Ks$ is antipodally
equivariant -- in particular if it is a connected 
sum of some fibered knot with\break%
itself -- then $\cone\pi$
can be taken to be a polynomial \
in $\xq$ alone. It was to this case\break%
that Looijenga drew explicit attention.) (2)~By 
replacing $\|\xq\|$ with a smooth,
monotone function of $\|\xq\|$ infinitely flat at 0, 
cf.\ \cite{Kauffman-Neumann}, $\cone\pi$ can be assumed\break%
smooth (but transcendental) at $\zero$.
\end{remarks}

\begin{proposition}\label{proposition 3.7} 
If $\Ks$ is a fibered link with smooth open-book structure $\pi$
then $\lambda(\Ks) = \lambda(\cone \pi; \zero)$ and 
$\rho(\Ks) = \rho(\cone \pi; \zero)$.  If
$\xq$ is an isolated critical point of $f$ and $f$ 
is tame at $\xq$ then $\lambda(\Ks(f;\xq)) = \lambda(f;\xq)$ 
and $\rho(\Ks(f;\xq)) = \rho(f;\xq)$.
\end{proposition}

\begin{proof}
We may assume $\pi$ is braided (\S\ref{section 2}). By 
taking the $\e$ in the definition of\break%
\newpage\noindent
``braided open-book structure'' sufficiently small, 
one may make the $2$-plane fields
$S_\Ks$ and $\Span{D(\cone\pi)}\mid S^3$ arbitrarily close; so 
they are homotopic.
\end{proof}

\smallskip

Let $\crit(f)$ denote the set of critical points 
of $f$. For each $P \in S^2$ consider the
sets
\begin{align*}
L^*(f, \pq) &= \crit(f)\cup (l\after \Span{Df})^{-1}(\pq), \\
R^*(f, \pq) &= \crit(f)\cup (r\after \Span{Df})^{-1}(\pq).
\end{align*}
If $A(\yq)$ and $B(\yq)$ denote the rows of $Df(\yq)$, 
considered as quaternions, then (cf.\ \ref{lemma 1.2})
\begin{align*}
L^*(f,\pq) \cup L^*(f,-\pq) &= 
    \{\yq: \Ps(\pq\Ps(B(\yq)(\conj A(\yq)))) = \zero\},\\
R^*(f,\pq) \cup R^*(f,-\pq) &= 
    \{\yq: \Ps(\pq\Ps((\conj A(\yq))B(\yq))) = \zero\},
\end{align*}
while 
$L^*(f,\pq) = \{\yq \in L^*(f,\pq) \cup L^*(f,-\pq):%
\pq\Ps(B(\yq)(\conj A(\yq))) \le 0$ (and so
on), when we identify $\Span{1}\sub\H$
with $\R$. Note that, for $\pq\ne \qq$, 
$L^*(f, \pq) \cap L^*(f,\qq) =\crit (f)=R^*(f,\pq)\cap R^*(f,\qq)$. 
(Note also that, though $L^*(f,\pq) \cup L^*(f,-\pq)$ 
and $R^*(f,\pq) \cup R^*(j,-\pq)$ are level sets 
of mappings to a $3$-dimensional
vectorspace, their expected codimension is not $3$ 
but $2$ because of the Pl\"ucker conditions.)

Now suppose $\crit(f)\cap D^4(\xq;\e)\sub \{\xq\}$.

\begin{hypothesis}\label{hypothesis 3.8}
$\pq,\qq \in S^2$, $\pq\ne\qq$, are such that 
(with respect to\break%
some convenient theory of 
geometric cycles representing ordinary homology
over $\Z$) 
the set $(l\after\Span{Df})^{-1}(\pq)\cap S^3(\xq;\e)$ 
(resp.\ $(l\after\Span{Df})^{-1}(\qq)\cap S^3(\xq;\e)$;
$(r\after\Span{Df})^{-1}(\pq)\cap S^3(\xq;\e)$;
$(r\after\Span{Df})^{-1}(\qq)\cap S^3(\xq;\e)$)
is the support of an absolute $1$-cycle in $S^3(\xq;\e)$ 
which bounds a relative $2$-cycle in $D^4(\xq;\e)$ supported by
$L^*(f,\pq)\cap D^4(\xq; \e)$ (resp.\ $L^*(f,\qq)\cap D^4(\xq;\e)$; 
$R^*(f,\pq)\cap D^4(\xq;\e)$; $R^*(f,\qq) \cap D^4(\xq;\e)$). 
(We will use the same symbols for the 
cycles and their supports.)
\end{hypothesis}

\begin{proposition}\label{proposition 3.9} 
Under Hypothesis \emph{\ref{hypothesis 3.8}}, $\lambda(f;\xq)$ 
\emph{(}resp.\ $\rho(f; \xq)$\emph{)} is the\break%
homological intersection number at $\xq$ of 
$L^*(f,\pq) \cap D^4(\xq;\e)$ and 
$L^*(f,\qq) \cap D^4(\xq;\e)$ 
\emph{(}resp.\ $R^*(f,\pq)\cap D^4(\xq;\e)$ and 
$R^*(f,\qq) \cap D^4(\xq; \e)$\emph{)}.
\end{proposition}

\medskip
\begin{proof}
This is a tautology, given the definitions 
of $\lambda$ and $\rho$ as Hopf invariants
and the relationship between linking numbers and 
intersection numbers.
\end{proof}
\medskip

\begin{remarks}\label{remarks 3.10} 
(1)~The point of \ref{hypothesis 3.8}-9
is that frequently \ref{hypothesis 3.8} can be verified, as,\break%
\newpage\noindent
for instance, in the examples in \S\ref{section 4}. 
(2)~One might conjecture that, for any $f$ with
an isolated critical point at $\xq$, \ref{hypothesis 3.8} 
holds for almost all pairs $(\pq,\qq)$. Certainly it
seems reasonable to expect, of a given $f$, that 
for almost all $\pq$ the sets $L^*(f, \pq)$
and $R^*(f,\pq)$ are ``2-manifolds with isolated 
singularities at $\xq$''. Perhaps some sort\break%
of higher-order tameness should be defined. 
(3)~If $f$ is a real-polynomial mapping,
then, for any $\pq$, $L^*(f,\pq) \cup L^*(f,-\pq)$ 
is a real-algebraic set and $L^*(f,\pq)$ is\break%
semi-algebraic (of course the same goes for $R^*$). 
Suppose $\xq$ is an isolated critical
point of $f$ and that $\pq,\qq \in S^2$, $\qq\ne\pm\pq$, are such 
that both $L^*(f,\pq) \cup L^*(f,-\pq)$ and\break%
$L^*(f, \qq)\cup L^*(f, -\qq)$ are purely $2$-dimensional 
near $\xq$. Then, near $\xq$, each of
$L^*(f, \pq)$, $L^*(f, -\pq)$, $L^*(f, \qq)$, and $L^*(f,-\qq)$ 
is the cone on some singuiar-link-\break%
with-integer-multiplicities, 
so \ref{hypothesis 3.8} holds. One might be tempted, therefore, to
reason from the distributive law that 
``$4\lambda(f,\xq)$ is the \emph{real-algebro-geometric}\break%
intersection number at $\xq$ of the real-algebraic 
surfaces $L^*(f, \pq) \cup  L^*(f, -\pq)$ and
$L^*(f,\qq) \cup L^*(f,-\qq)$''. It seems hard to make that 
statement true inside real algebraic geometry! 
(Real-algebraic cycles are naturally oriented over $\Z/2\Z$\break%
rather than over $\Z$. In the present case, even if 
each of the algebraic sets\break%
$L^*(f, \pq) \cup L^*(f, -\pq)$ 
is purely $2$-dimensional, giving them local $\Z$-orientations
algebro-geometrically is complicated by the fact that the parameter 
space for this family of surfaces is the non-orient able real 
projective plane $RP^2$ rather than $S^2$;\break%
cf.\ \hyperlink{twist}{the last sentence of \ref{machinery 4.1}}.) 
Perhaps there exists (I have not been able to learn of
it) an applicable theory of integer intersection 
numbers, \emph{calculable inside real\break%
semi-algebraic geometry}, and giving the correct 
topological answers? (4)~As\break%
mentioned in \S\ref{section 0}, 1 know (as of December, 1986) of 
no example of a fibered link
$\Ks$ with $\lambda(\Ks) < 0$. Especially if no such link 
exists, it would be interesting to
known whether there exists a function $f$ with an 
isolated critical point $\xq$ such that
$\lambda(f; \xq) < 0$.\hyperlink{regret}{$^{(*)}$}
\end{remarks}

\medskip

\section{Examples}\label{section 4}

Most of the examples in this section involve 
complex analyticity somehow, so
we begin by introducing some complex machinery.

\begin{machinery}\label{machinery 4.1}
At a point where $F:\C^2\to\C$ is smooth, the 
\emph{complex differential} $D_\C F$ is the complex row vector 
$[F_z\, F_{\bar z}\, F_w\, F_{\bar w}]$, where 
$F_z = (F_x - iF_y)/2$,
$F_{\bar z} = (F_x + iF_y)/2i$, etc., and subscripts 
indicate partial differentiation. In terms of\break%
$D_\C F$, the real differential matrix $DF$ is
\[
\begin{bmatrix}
\Re(F_z+F_{\bar z}) & \Re(iF_z-iF_{\bar z}) &
               \Re(F_w+F_{\bar w}) & \Re(iF_w-iF_{\bar w}) \\
\Im(F_z+F_{\bar z}) & \Im(iF_z-iF_{\bar z}) &
               \Im(F_w+F_{\bar w}) & \Im(iF_w-iF_{\bar w}) 
\end{bmatrix}.
\]

\newpage\noindent
As in \ref{lemma 1.2}, we see that at a regular point of $F$, 
$l(\Span{DF})$ is the unit vector of
\begin{equation}\tag{$*$}\label{equation star}
(|F_z|^2 - |F_{\bar z}|^2 + |F_w|^2 - |F_{\bar w}|^2)\iq %
-2\Im(F_z \overline{F_{\bar w}} - \overline{F_{\bar z}}F_w)\jq %
-2\Re(F_z \overline{F_{\bar w}} - \overline{F_{\bar z}}F_w)\kq; 
\end{equation}
similarly, at a regular point of F, 
$r(\Span{DF})$ is the unit vector of
\begin{equation}\tag{$**$}\label{equation starstar}
(|F_z|^2 - |F_{\bar z}|^2 + |F_{\bar w}|^2 - |F_w|^2)\iq %
-2\Im(F_z\overline{F_w} - \overline{F_{\bar z}}F_{\bar w})\jq %
-2\Re(F_z\overline{F_w} - \overline{F_{\bar z}}F_{\bar w})\kq; 
\end{equation}
and (if $F$ is smooth everywhere) $\crit(F)$ is 
defined by the vanishing of either (\ref{equation star})
or (\ref{equation starstar}).

In this complex context, we will have a 
particular interest in $L^*(F,\iq) \cup L^*(F, -\iq)$ 
and $R^*(F,\iq) \cup  R^*(F,-\iq)$. As sets,
\begin{equation*}
L^*(F,\pm\iq) = \{F_z \overline{F_{\bar w}} - \overline{F_{\bar z}}F_w =0,
\pm(|F_z|^2 - |F_{\bar z}|^2 + |F_w|^2 - |F_{\bar w}|^2)\ge 0\}
\end{equation*}
and
\begin{equation*}
R^*(F,\pm\iq) = \{F_z \overline{F_{w}} - \overline{F_{\bar z}}F_{\bar w} =0,
\pm(|F_z|^2 - |F_{\bar z}|^2 - |F_w|^2 + |F_{\bar w}|^2)\ge 0\}.
\end{equation*}
Suppose $F_z \overline{F_{\bar w}} - \overline{F_{\bar z}}F_w$
and $F_z \overline{F_{w}} - \overline{F_{\bar z}}F_{\bar w}$
are products of complex analytic functions
and conjugates of complex analytic functions. 
Then their level sets, where they
are $2$-dimensional, are equipped with natural 
integer multiplicities; in particular
this is true of the sets of zeroes, and so at any 
isolated critical point of $F$ near
which $L^*(F,\pm\iq)$ and $R^*(F,\pm\iq)$ are $2$-dimensional, 
Hypothesis \ref{hypothesis 3.8} will be satisfied
(with $\pq = \iq$, $\qq = -\iq$). Note, however, that 
\hypertarget{twist}{the multiplicity assigned by the defining
function must be twisted} by the sign of $\iq$ to give 
the multiplicity needed for \ref{hypothesis 3.8}
(consider the local coordinates on $S^2$ given by 
stereographic projection from the two poles $\iq$ and $-\iq$).
\end{machinery}

\begin{example}\label{example 4.2} Let $f\from \C^2\to\C$
be a complex polynomial. If $f$ is squarefree,
then any critical point $(z, w)$ is necessarily 
isolated. Claim: in this case,\break%
$\lambda(\Ks(f;(z, w))) = 0$. Proof: at any regular point 
of $f$, $\ker Df$ is a complex line, so
$l\after\Span{Df}$ is identically $\iq$ 
and $\lambda(\Ks(f;(z, w)) = \lambda(f;(z, w))$ 
is the Hopf invariant of\break%
a constant.

It follows from \ref{corollary 2.6} that 
$\rho(\Ks(f;(z, w)))=\mu(\Ks(f;(z, w)))$. 
In fact, $R^*(f; -\iq)$\break%
is the complex plane curve 
$\{f_z=0\}$ with the opposite orientation to that given by
its complex structure, and $R^*(f;\iq)$ is the conjugate-complex 
plane curve $\{\bar{f_w}=0\}$ with the orientation given by its 
conjugate-complex structure; the intersection
number at $(z, w)$ of these cycles is then $(-1)\cdot(-1) = 1$ 
times the intersection\break%
\newpage\noindent%
number at $(z, w)$ of the complex plane curves 
$\{f_z=0\}$ and $\{f_w=0\}$, which is
Milnor's definition of $\mu(\Ks(f;(z, w)))$ 
\cite[p.\ 59]{Milnor:singular-points}.
\end{example}

\begin{remarks}\label{remarks 4.3} (1)~The links $\Ks(f;(z, w))$ 
are well understood (cf. \cite{Eisenbud-Neumann}, 
\cite{Le:sur-les-noeuds},\break%
etc.); they are (quite restricted) iterated torus 
links, and also (very special) closed
strictly positive braids. In \cite{Rudolph:isocp2} it is shown 
that $\lambda(\Ks) = 0$, $\rho(\Ks) = \mu(\Ks)$ for any\break%
closed strictly positive braid $\Ks$. (More 
generally, if $\Ks$ is a closed strictly\break%
homogeneous braid \cite{Stallings} or even a closed 
generalized strictly homogeneous braid,
then $\lambda(\Ks)$ and $\rho(\Ks)$ are the negative and 
positive parts of $\mu(\Ks)$ in an obvious
sense.) In \cite{Neumann-Rudolph:computing}, 
$\lambda$ and $\rho$ are calculated for all 
fibered iterated torus links. (2)~\ref{example 4.2}
substantiates the intuition that the link of a 
complex plane curve singularity (or
any closed strictly positive braid) is somehow 
``as positive as it can be''. It should
be contrasted with the fact that, though the 
symmetrized Seifert form of such a
link has non-negative signature, \cite{Rudolph:signature}, it is only 
rarely positive-definite -- for\break%
complex plane curves, this happens exactly when 
the singularity is ``simple'' in the
sense of Arnol'd. (Actually, the sign convention 
in \cite{Rudolph:signature} is unusual; with the more
standard one, a closed positive braid has non-positive signature.)
\end{remarks}

\begin{example}\label{example 4.4} 
Let $\Rev\from S^3\to  S^3$ be an orientation-reversing 
diffeo-\break%
morphism. The mirror image of a link $\Ks = (S^3, K)$
is the link $\Rev \Ks =\break%
(S^3, \Rev K)$. 
Claim: if $\Ks$ is fibered, then 
$\lambda(\Rev \Ks) = \rho(\Ks)$ (so also 
$\rho(\Rev \Ks) = \lambda(\Ks))$. 
Proof: this is a simple calculation from 
the formulas in \ref{theorem 2.3} and 
\ref{theorem 2.5}-6.\break%
(More generally, if $f\from \H\to\H$ has an isolated 
critical point at $\zero$, then $\conj\after f$ has\break%
also and $\lambda(\conj\after f; \zero) = \rho(f; 0)$, by 
consideration of the effect of $\conj$ on $\pi_3(G)$;
of course, $\Ks(\conj\after f; \zero)$ is a mirror image of 
$\Ks(f;\zero)$.) In particular, if $\Ks$ is
amphicheiral (i.e., isotopic to its mirror image) 
then $\lambda(\Ks) = \rho(\Ks)$.
\end{example}

\begin{example}\label{example 4.5} 
The figure-$8$ knot $\Ks$ is amphicheiral, 
and $\mu(\Ks) = 2$, so by\break%
\ref{example 4.4}, $\lambda(\Ks) = 1 = \rho(\Ks)$. 
Now, $\Ks$ is a closed homogeneous braid (the closure of the
homogeneous braid word $\s_1\s_2^{-1}\s_1\s_2^{-1}$ in the $3$-string 
braid group $B_3$), and the
techniques of \cite{Rudolph:isocp2} could also be brought to bear. 
But it is most entertaining to
calculate $\lambda(\Ks)$ and $\rho(\Ks)$ by \ref{machinery 4.1}.

First let $F(z,w) = w^3 - 3|z|^2(1 + z -\bar{z})w - 2(z + \bar{z})$. 
Then $D_\C F(z, w) = %
\begin{bmatrix} -3w(\bar{z}+2|z|^2 - {\bar{z}}^2)-2
              & -3w(z + z^2-2|z|^2)- 2 
              & 3w^2-3|z|^2(1+z-\bar{z}) \medspace 0
\end{bmatrix}$ so\break%
$L^*(F,\iq) \cup L^*(F,-\iq) = %
\left\{\left(-3\bar{w}(\bar{z} + \bar{z}^2 - 2|z|^2)-2\right)%
\left(w^2 - |z|^2(1 + z -\bar{z})\right) = 0\right\}$
and
$R^*(F,\iq) \cup R^*(F,-\iq) = %
\left\{\left(-3{w}(\bar{z} + 2|z|^2 - \bar{z}^2)-2\right)%
\left(\bar{w}^2 - |z|^2(1 + z -\bar{z})\right) = 0\right\}$. For
$|z|^2 + |w|^2$ small, 
$\left(-3\bar{w}(\bar{z} + \bar{z}^2 - 2|z|^2)-2\right)%
\left(-3{w}(\bar{z} + 2|z|^2 - \bar{z}^2)-2\right)\ne 0$,
so at a\break%
point of $\crit(F)$ near $(0,0)$, 
$w^2 = |z|^2(1+z-\bar{z})$, $w =\pm|z|+o(|z|^3)$; then
$|F_z|^2 - |F_{\bar z}|^2 + |F_w|^2 - |F_{\bar w}|^2 %
=\pm 4 |z| \Re (\bar{z}-z + 4|z|^2-\bar{z}^2-z^2)+o(|z|^4)=$\break%
\newpage\noindent
$\pm 8 |z| \left((\Re z)^2 + 3(\Im z)^2\right) + o(|z|^4)$.  Since
$\left((\Re z)^2 + 3(\Im z)^2\right)$ is positive-definite, 
the critical point of $F$ at $(0, 0)$ is isolated. 
Also, $L^*(F,\pm\iq)$ (resp.\ $R^*(F,\pm\iq)$) is\break%
well approximated near $(0,0)$ by $\{w = \pm|z|\}$
(resp.\ $\{w =\mp |z|\}$). These cycles
have intersection number $0$ at $(0,0)$, so 
$\lambda(F; (0,0)) = 0 = \rho(F; (0, 0))$, so\break%
$\mu(\Ks(F;(0, 0))) = 0$; though $(0, 0)$ is a 
genuine critical point, $\Ks(F; (0, 0))$ is\break%
unknotted.

Let $G(z,w) = F(z^2, w)$; again $(0,0)$ is an 
isolated critical point; 
\begin{align*}
L^*(G, \iq)\cup L^*(G, -\iq) = & 
      \{(2\bar{z} [-3\bar{w}(\bar{z}^2 + \bar{z}^4 - 2|z|^4)- 2]) \\
    & \times (w^2-|z|^4(1 + z^2 - \bar{z}^2)) = 0\}
\end{align*}
and
\[
R^*(G, \iq)\cup R^*(G, -\iq) = %
    \{(2{z} [-3{w}(\bar{z}^2 + 2|z|^4 - \bar{z}^4)- 2])%
      (w^2-|z|^4(1 + z^2 - \bar{z}^2)) = 0\}
\]
so one quickly calculates 
$\lambda(G; (0, 0)) = 1 = \rho(G; (0,0)$, 
$\mu(\Ks(G; (0,0/)) = 2$.
Now, $K(G; (0, 0))$ double-covers $K(F;(0,0))$, which 
is connected, so it has $1$ or $2$\break%
components -- but 
$\mu(\Ks(G; (0,0))$ is even, so $K(G; (0,0))$ has an 
odd number of components. Thus $K(G; (0, 0))$ is connected and 
$\Ks(G; (0, 0))$ is a knot. It must\break%
be the figure-$8$ knot. (Only three fibered knots 
have Milnor number $2$ -- the two trefoils and the 
figure-$8$ knot. One trefoil is $\Ks(z^2 + w^3; (0, 0))$; 
$\lambda(\Ks(z^2 + w^3; (0, 0)) = 0$ by \ref{example 4.2}. 
The other trefoil is $\Rev \Ks(z^2 + w^3; (0, 0))$; 
$\rho(\Rev \Ks(z^2 + w^3; (0, 0)) = 0$ by \ref{example 4.4}.)

Of course it is easy enough to see directly that 
$\Ks(G; (0, 0))$ is a figure-$8$ knot,
by considering the closed braid cut out by 
$G=0$ in a sufficiently small bidisk
boundary $\{(z, w): |z|\le\e, |w|\le \e'\}$,
which is readily seen to be the closure of
$\s_1\s_2^{-1}\s_1\s_2^{-1}$.
\end{example}

\begin{remark}\label{remark 4.6} 
Perron was the first to give a 
real-polynomial mapping $\R^4\to\R^2$\break%
having an isolated critical point with 
local link the figure-$8$ knot, \cite{Perron:figure8}. His\break%
polynomial is somewhat more complicated than that 
in \ref{example 4.5}, and in particular has
resisted my occasional attempts to use it to 
calculate $\lambda$ and $\rho$; the ``half-\break%
complex'' nature of $F$ and $G$ (the vanishing of 
their $\bar{w}$-derivatives) is a great\break%
simplification.
\end{remark}

\begin{example}\label{example 4.7} 
Let $f(z, w) = z^2 + w^3$, $g(z, w) = z^3+ w^2$, $F = fg$. 
Then
\begin{align*}
& D_\C F(z,w) = %
\begin{bmatrix} 2z\bar{g} & 3\bar{z}^2f %
                    & 3w^2\bar{g} & 2\bar{w}f\end{bmatrix},\\
& L^*(F,\iq)\cup L^*(F,-\iq)=\left\{(4zw-9z^2 w^2)\bar{f}\bar{g}=0\right\},\\
& R^*(F,\iq)\cup R^*(F,-\iq)=\left\{6z\bar{w}^2|g|^2-6z^2\bar{w}|f|^2=0\right\}.
\end{align*}
The origin is an isolated critical point of $F$. 
\ref{machinery 4.1} applies directly to calculate\break%
$\lambda$: $L^*(\iq) = \{z\bar{f}=0\}$ and 
$-L^*(F,-\iq) = \{w\bar{g}(4-9zw)=0\}$ as cycles, so\break%
$\lambda(F; (0, 0)) = 1$. (Of course the complex curve 
$4-9zw = 0$ doesn't pass through $(0,0)$ so it isn't 
involved in the calculation.) \ref{machinery 4.1} 
doesn't quite apply to calculate $\rho$, 
since (as a set) $R^*(F, \iq) = \{\bar{w}|g|^2-z|f|^2=0\}$,
and $\bar{w}|g|^2-z|f|^2=0$ isn't just the
product of some complex analytic and some 
conjugate-analytic factors. But one
may verify that there is a neighborhood of $(0, 0)$
which has the same intersection
with $R^*(F, \iq)$ as it has with $w -\bar{z} = 0$, and 
then calculate $\rho(F; (0, 0)) = 2$ (since
$-R^*(F,-\iq) = \{z\bar{w} = 0\}$ as a cycle).

This example can be generalized. 
Let $a, b, c, d$ be positive integers,\break%
$f(z, w) = z^a + w^b$, $g(z,w) = z^c + w^d$, 
$F(z, w) = f(z, w)g(\bar{z},\bar{w})$, 
$G(z, w) =$\break%
$f(z, w)g({z},\bar{w})$.  Then 
$\Ks(F;(0,0))$ and $\Ks(G;(0,0))$ are certain iterated torus\break%
links (of $\GCD(a,b) + \GCD (c, d) > 1$ components). 
For most choices of $a$, $b$, $c$,\break%
$d$, the critical point of 
$F$ (resp.\ $G$) at $(0, 0)$ is isolated so $\Ks(F; (0, 0))$
(resp.\ $\Ks(G; (0,0))$) is a fibered link; the invariants 
$\lambda(F; (0,0)), \dots, \rho(G; (0,0))$ can be
calculated. Typically such a link is neither 
(isotopic to) the link of a complex\break%
plane curve singularity nor (isotopic to) the 
mirror image of such a link; this is\break%
detected by $\lambda$ and $\rho$ without recourse to the 
classification of links of curve\break%
singularities. Note that for certain bad choices 
of exponents, the critical point of $F$\break%
or of $G$ at (0,0) will not be isolated; e.g.,
$a = c$ is bad for $F$, and $b = d$ is bad for\break%
both $F$ and $G$. Note also that it can be determined 
just which of the links $\Ks(F; (0,0))$ and $\Ks(G; (0,0))$ are, 
and are not, fibered -- for instance, by using\break%
the calculus of splice diagrams \cite{Eisenbud-Neumann}. 
Interestingly, it appears that whenever\break%
$\Ks(F; (0, 0))$ is fibered, in fact $F$ has an 
isolated critical point at $(0,0)$, and\break%
likewise for $G$ (cf.\ \cite{Rudolph:mbob}).
\end{example}
\bigskip

\renewcommand\S{\relax}

\section*{Acknowledgements}

This research was supported by the Fonds National 
Suisse and by MSRI. I am\break%
happy to acknowledge the support of my friends at 
the University of Geneva,
especially Claude Weber and Gerhard Wanner, 
quondam Pr\'esident de la Section
de Math\'ematiques. At MSRI I had helpful 
conversations with many colleagues,
particularly Walter Neumann.
\bigskip

\makeatletter
\renewcommand\section{\@startsection{section}{1}%
  \z@{1.7\linespacing\@plus\linespacing}{1.5\linespacing}%
  {\normalfont\mdseries\scshape\noindent}}
\makeatother

\renewcommand{\bysame}{\leavevmode\raise.5ex\hbox to2em{\hrulefill}%
\thinspace}

\smallskip

{\footnotesize
{\obeylines
{\emph{
\noindent Box 251, Adamsville
\noindent Rhode Island 02801 U.S.A.
\noindent ~
\noindent Current address:
\noindent Department of Mathematics
\noindent Clark University, Worcester
\noindent Massachusetts 01610 U.S.A.}
}
\smallskip
\noindent Received August 27, 1985/January 20, 1987
}
}
\smallskip

\begin{addenda}\label{addenda}
\hypertarget{projected}{(1)}~The 
\hyperlink{projected-backlink}{extension to non-fibered links}
remains to be done. ~\hypertarget{more Gabai}{(2)~Since published}, 
\cite{Gabai:more}.~\hypertarget{regret}{(3)}~In fact,
\hyperlink{regret-backlink}{many links have $\lambda(\Ks)<0$}; 
an explicit example (actually from the 
family of iterated torus links studied in \ref{example 4.7}) 
is worked out in \cite{Neumann-Rudolph:diff-index}.
~\hypertarget{more N-R}{(4)~Since published}, 
\cite{Neumann-Rudolph:enhance-higher}--\cite{Neumann-Rudolph:diff-index}.
\end{addenda}

\renewcommand{\refname}{\protect\vskip-.5in\relax}

\end{document}